\documentclass[czech]{article}

\begin{document}

\title{The trace of First World War on mathematics in Brno}
\author{ Laurent MAZLIAK\footnote{Laboratoire de Probabilit\'es et Mod\`eles al\'eatoires \& Institut de Math\'ematiques (Histoire des Sciences Math\'ematiques), {Universit\'e Paris  VI,
France}. mazliak@ccr.jussieu.fr}
 and  Pavel \v{S}I\v{S}MA\footnote{Katedra matematiky P\v{r}\'\i rodov\v{e}deck\'e fakulty Masarykovy univerzity, { Brno,  Czech Republic}. sisma@math.muni.cz}
}

\maketitle 

\begin{quote} {\bf Abstract} : In the present paper, we study the evolution of the mathematical community in Brno (Br\" unn in German), the Moravian maintown, in the years between 1900 and 1930. In particular, we want to discuss how the Great War and its consequences (creation of Czechoslovakia and shift of power from the Germans to the Czechs) had an effect on this evolution.



\end{quote}

\section*{Introduction}
The foundation of Czechoslovakia in 1918 appears a good example of an attempt of reorganization of Europe after the end 
of World War I. In places where there was a tense cohabitation of several national communities (as in many parts 
of the collapsed Austro-Hungarian empire), it was necessary to choose a form of organization allowing the coexistence 
of several traditions. This was in particular the case with the organization of educational system. 

 Brno (Br\"unn in German)\footnote{To simplify, we shall 
always use the Czech name of the city when referring to it, though, when the context indicates that the German community 
is concerned, it obviously must be replaced by Br\" unn.}, the capital of the Moravia district was such a place of cohabitation and appears therefore as a good laboratory to understand the kind of 
ruptures and continuities in history around WW1. Moravia is a border region of 
Austria, and therefore has always been a crossroads of cultures. Due to the presence of a very strong 
German minority, Moravia was in 1918 one of the parts of Europe where the question of nationalities would appear with particular 
acuteness. The father of Czechoslovakian independence, Tom\'a\v{s} G.Masaryk wrote in {\it New Europe} that in the 
{\it so-called `German territories' in Bohemia (Moravia and Silesia) live numerous Czechs; it is therefore fair that the 
Czech state will keep them; it would be unfair to sacrifice hundreds thousand Czechs to the 
{\rm furor teutonicus}}\footnote{Masaryk, T. G. : Nov\'a Evropa. Stanovisko slovansk\'e. (New Europe. A slavonic point of view), 
Dubsk\'y, Praha, 1920. 
}.
This peremptory assertion hardly dissimulated that one may expect serious problems from this ambiguous annexation, 
an impression confirmed by Bene\v{s}'s declarations. Bene\v{s}, the minister of foreign affairs of the new Czechoslovakia, came to 
Paris on 5 February 1919 in order to explain the proposed drawing of borders of the new state to the delegates 
preparing the Peace conference. He mentioned that {\it the relations of Czechoslovakia with its neighbours have to be 
settled in order to avoid any conflict}. To achieve this goal Bene\v{s} found {\it necessary to take the ethnographic map 
into consideration with  maximal care, above all in regions where natural borders are not obvious}\footnote{Quoted in the 
newspaper Le Matin, 6 February 1919.}. 

To understand the case of Brno, it is therefore
vital to understand how the difficult contacts between the Czech majority and the large German minority had influenced 
the whole process of edification of the education institutions between ca 1880 and 1930. Though the German minority  
lost its domination in 1918, the institutions were still much influenced by the culture that had prevailed before the war, 
though there were several attempts to create a new interest towards the countries of the victorious side. 

And hence an intricate question (maybe the most difficult one of the present paper) immediately arises in what we are now 
writing. Namely : who are, in this story,  the Germans, who are the Czechs? What makes the question hard is that the 
answer was not univocal and was changing during the period depending on political and social conditions. K\v ren has 
devoted a study\footnote{K\v ren, J. : Changes in identity. Germans in Bohemia and Moravia in the nineteenth and twentieth 
centuries. In Bohemia in History (Mikul\'a\v s Teich,  ed.), Cambridge University Press, 1998.}  to the special case of 
the `Germans' which reveals how the definition of who was German (and therefore who was Czech) fluctuated. In the population censuses of the years 1880-1900,  the numbers of inhabitants of Brno who declared to have Czech as their 
language of communication, varied from 30 to 40 percent. In the last census before WW1 (1910), 41.000  out of 126.000  
Brno inhabitants (32\%) declared Czech to be their usual language\footnote{D\v r\'\i{}mal, J., Pe\v sa, V. : D\v ejiny m\v esta Brna, vol. 2. Blok, Brno
1973. p. 64.}.  These figures must be considered with care.  Political and economical reasons probably lead to an overestimation of the 
German settlement because in January 1919, 61\% of the (almost identical) inhabitants declared to belong 
to the Czech community. 

Our study emphasizes the fact that may seem obvious at first glance. The history of mathematics, and more widely the 
history of intellectual life in Moravia cannot avoid taking into account the question of  relationship between the 
two national communities of German inhabitants on the one hand, and Czech inhabitants on the other hand, even if the 
contours of these communities were never very precise. These  complicated 
contacts, mixing rivalry (often) and dialogue (sometimes) appear as a basic riddle by means of which one can explain much of 
the history of academic life in the Czech lands during two centuries.  To understand  the situation   
of 1918, it is, however, necessary to consider this question with care and to overcome the impressions created by our 
knowledge of the violent end of this cohabitation. The German invasion of 1939 followed by the terrible years of 
occupation,  and the general expulsion of German inhabitants between 1946 and 1948 
might lead us to imagine the presence of the two communities as a permanent battlefield. In fact, it rather seems that 
 the members of both communities had mostly ignored each other before 1918 - and hence passively accepted living 
in the Habsburg monarchy for granted.  As Fejt\"o mentions, {\it it would be projecting oneself abusively in the past to believe that, even at this moment [in 1916, at the death of emperor Franz-Josef], the 
majority of Czechs were ready to get rid of the monarchy}\footnote{Fejt\"o, F. : Requiem pour un Empire d\'efunt. 
Histoire de la destruction de l'Autriche-Hongrie, Points-Histoire, Seuil, 1993. p.141.}. We, therefore, think that 
the mentioned riddle is appropriate to describe not only the splits, but also the continuities in academic life.

\section{The Czech fight for higher education in Moravia before WW1}

The first real university in Brno,  a Technical University ({\it Technische Hochschule}) was founded in Brno in 1873, replacing the Polytechnicum established in 1849, itself a distant heir of the old Olomouc Academy. The Technical University 
was divided into `faculties' and was managed by an elected rector. Though the number of
professors increased, the number of students stagnated and the Brno Technical
University was, in fact, a small institution\footnote{Hellmer, K. : Geschichte der Deutschen Technischen
Hochschule in Br\"unn. In Festschrift der k. k. Technischen Hochschule in Br\"unn
zur Feier ihres f\"unfzigj\"ahrigen Bestehens und der Vollendung des
Erweiterungsbaues im October 1899. Br\"unn 1899.}. Numerous Austrian technicians and scientists began their academic 
career in modest size institutions. Havr\'anek\footnote{Havr\'anek, J. : The university professors and students in nineteenth-century Bohemia. In Bohemia in History (Mikul\'a\v s Teich,  ed.), Cambridge University Press, 1998.} quotes the epigrammatic characterization of the 
professor's career in the Habsburg monarchy : {\it Sentenced to Czernowitz, pardoned to Graz, promoted to Vienna}. Though 
less prestigious than Graz because of the `hostile' Czech environment, the Brno Technical University was certainly seen as a reasonable position because of its proximity to Vienna.  Such a position was considered  
a springboard for accessing one of the Vienna universities. 

Since the beginning of the Czech national revival (i.~e. the end of the 18th century), the `war' between the communities had been waged mainly 
on a symbolic ground, namely in cultural and intellectual life. The intellectual conflict between the Germans and Czechs 
existed obviously in every place where the two communities cohabited. It is worth noticing that there was a noticeable difference between the situation in Bohemia and in Moravia. 

In Prague, due to the large numerical domination of the Czech community, the Germans often live the cohabitation as a threat. The division of the Prague university in 1882  into German and  Czech Universities was seen by the Germans as a necessity to preserve German high education in Prague. The conflicts in Prague were often violent. An example is the creation of the movement {\it Los von Prag} (Free or Out from Prague), where Germans from northern and western Bohemia, resenting Prague's anachronistic liberalism, attempted to move the cultural and ideological centers of German Bohemia to Reichenberg (Liberec). In 1897 German students attempted to move German universities from Prague to Liberec\footnote{See Gary B. Cohen: The Politics of Ethnic Survival: Germans in Prague, 1861-1914, Purdue University Press, 2006.}. Many testimonies exaggerating these tensions at 
the beginning of 20th century can be found  in the intellectual life. A famous example is given by the celebrated 
Ha\v sek's book {\it The good soldier \'Svejk}\footnote{Ha\v{s}ek, J. : The Good Soldier Schweik, tr. Paul Selver, Boni, 
New-York, 1930.} where the Austrian militarism and bureaucracy is ridiculed. On the academic ground, we find a harsh 
dialogue between two personalities of Prague intellectual life, the rector of the (German) Charles University A.Sauer  
and the politically involved Czech physician O.Srd\'\i{}nko. In 1907, mocking the complaints of the Prague Czech University 
members to obtain better financial subsidies, the rector Sauer called the latter institution a {\it spoiled child}, 
bringing Srd\'\i{}nko's wrath on him in a small brochure published in 1908\footnote{Srd\'\i{}nko, O. :
Zh\'y\v ckan\'a \v cesk\'a universita. Odpov\v e\v d prof. Dr. A. Sauerovi, rektoru (1907--8)
n\v emeck\'e university v Praze. (The spoiled Czech University. 
An answer to Prof. A.Sauer, rector of the German University in Prague), N\'akladem vlastn\'\i{}m, Praha, 1908.}. 
 
 The already 
mentioned geographic situation of south-Moravia on the border of cultures and languages created a different mode of cohabitation. In Brno,  the Germans, though not in the majority, formed a strong community. So close to the imperial maintown, they were in position to keep a tight control on the evolution of the situation. The division of the Prague University had resulted in an increased capacity 
at both universities. The Czechs from Moravia went to the Czech Universities in Prague; the Germans mostly studied  in
Vienna. At the end of the 1880s,  about 1000 students from Moravia studied in Austrian universities (700 in Vienna, 250 at 
the Czech University in Prague and 60 at the German University in Prague). In the middle of the 1890s, the number of 
students from Moravia increased to 1300. The 
question of creating a new university in Moravia, therefore, reappeared. Moreover, the existence of tiny universities in the Habsburg monarchy in places where the students potential was much smaller, such as  the university in Czernowitz in distant Bukovina which had less than 300 students, was seen as a proof for the  soundness of the project\footnote{Jord\'an, F. : D\v ejiny
university v Brn\v e. Universita J. E. Purkyn\v e v Brn\v e. Brno 1969. pp. 40, 43.}. 

The central problem was the teaching language. 5 million Czechs had only one university in Prague, whereas 8 million 
Germans had 5 Universities (Vienna, Prague, Graz, Innsbruck, and Czernowitz). The Czechs from Moravia and Silesia formed 
a large group among the university students in Austria.  In 1885,  Masaryk mentioned the necessity for 
the Czech nation to require two universities for the first time\footnote{Masaryk, T. G. : Jak zvelebovati
na\v si literaturu naukovou, Athenaeum II, 1885, p. 275.}.  
He wrote later that this was necessary in order to guarantee a wealthy competition and new places for young 
teachers\footnote{Masaryk, T. G. : Druh\'a universita, Na\v se doba, I, 1894,
pp. 672--676.}.

Masaryk brought the following reasoning to support his request. If a young specialist were appointed to the only  
university of the Czech lands, in the next decades, no other talented young man had any  chance to obtain a professorship 
at university. As a confirmation, at the end of 1890s, a high number of privatdocents at the Czech University in Prague 
resulted in an academic jamming to be solved.  The Vienna government adopted a wait-and-see policy. Already in 1896, 
the minister Gautsch conceded a right of Moravia for a Czech university and a Czech technical university, but it 
remained a formal declaration as Vienna asked for a preliminary agreement about the 
creation of a new university from the national components\footnote{Jord\'an, F. : D\v ejiny
university v Brn\v e. Universita J. E. Purkyn\v e v Brn\v e. Brno 1969. pp. 47--48.}. 

 In 1899, a Czech Technical University was opened in Brno as a counterpart of the German Technical University but the foundation of a Czech university was opposed to. The German politicians  and intelligentsia in the 
Czech lands warned against its establishing. 
They emphasized that a Czech university would not be in position to prepare the first rate specialists because the Czech 
language limited contacts with the world leading scientists\footnote{Bachmann,
A. : Die Universit\"aten im \"osterreichischen V\"olkerstreite. Neue Freie Presse,
1st November 1902.}. But they also expressed their fear of a forced `czechisation', 
sometimes in harsh terms. A. Sauer wrote in 1907 : {\it And the students of the 
university of Br\"unn, objects of every wish, where are they going to find employment, if it is not by occupying those 
that, until today, have been in the hands of representatives of other races [sic] , and above all of the Germans? [\dots ] 
If czechisation means expelling the Germans from their positions, dominating the Germans, oppressing the Germans, then 
this proud slogan must not ornate the pediment of a second Czech university without our deepest protesting against this 
villainy}\footnote{ A.Sauer. Quoted by O.Srd\'\i{}nko, op. cit.}. Moreover, Sauer mentioned that a major part of the taxes were 
paid by the (richer) German community who would be therefore reluctant to finance a second Czech university. This remark 
made Srd\'\i{}nko explode :   {\it He who wants to make an object of trade with civilization and to offer it only to the rich, 
behaves himself as a barbarian, even if he is endowed with the dignity of a university rector}\footnote{Srd\'\i{}nko, O., op.cit.}.

In 1905, Masaryk talked again about the necessity of a second university for improving the quality of the first 
University. He wrote moreover that the second Czech university had to be in Brno\footnote{Lidov\'e
Noviny, 28 June 1905.}. The year 1905 was the top of the 
efforts to establish a Czech university in Brno. Yet, after 1908 no Czech political party kept  it in its political 
program, though Masaryk presented  a series of petitions\footnote{Jord\'an, F. : D\v ejiny university v Brn\v e. Universita
J. E. Purkyn\v e v Brn\v e. Brno 1969. p. 105.} 
 for the creation of a new Czech university at the Austrian Parliament in June 1912. All facts considered, due to the political 
situation in Brno and the tergiversations of Vienna government, the attempts for its creation resulted in a failure. 

\section{Mathematics in Brno at the eve of WW1}

Let us now draw the picture of mathematics in higher education in Brno at the eve of WW1 (between 1880 and 1914).  The comparison with the picture of the post-war situation described in the subsequent 
section should emphasize the continuities and the splits in the local mathematics between the two periods. We shall successively consider the case of the German and of the Czech institutions.

\subsection{German institutions}

At the German Technical University,  the chairs ({\it Lehrstuhls}) of mathematics,  of descriptive geometry, 
and of theoretical mechanics met a more frequent the turn-over than in the other departments. This may have been the case 
mainly  because mathematicians had more opportunities to find another position  (mathematics was taught both in 
technical universities and in universities). Also, these mathematicians mostly came to Brno from Graz, Prague and 
Vienna where they were born, had relatives and studied at local universities. Many of them, therefore, used the first opportunity 
to come back\footnote{\v Si\v sma, P. :
Matematika na n\v emeck\'e technice v Brn\v e. Prometheus, Praha 2002.}.

The first renowned  mathematician teaching at the Brno Technical University
 Emanuel Czuber (1851--1925)\footnote{Dolezal, E. : Emanuel Czuber.
Jahresbericht der Deutschen Mathematiker--Vereinigung. 37, 1928, pp.
287--297.}, who came from the Prague German Technical University in 1886, is a good example of how Brno was used as the aforementioned springboard to Vienna. Czuber was a Czech from Prague (his original name,  \v{C}ubr, was germanized during his studies at the German Technical University in Prague).   In 1891, he left Brno for the Vienna Technical University. 

In 1906, a brilliant period for mathematics opened at the German Technical University in Brno, which  
was probably directly linked with the position obtained by Georg Hamel (1877--1954) in 1905 after Hellmer's retirement. This was the 
first post the 28 year-old Hamel received after his doctorate thesis under Hilbert in G\"ottingen (about Hilbert's 
fourth problem) and his habilitation in Karlsruhe. Hamel's presence and unprecedented activity (in particular for 
recruiting the first rate collaborators) seems to have boosted the mathematical life in Brno from 1905 as is 
testified by the sudden appearance of a Brno seminar in the journal of the German mathematical society ({\it Deutsche 
Mathematiker--Vereinigung, DMV}). In the Czech lands at the end of the 19th century two mathematical unions coexisted. 
There was the German Union of Mathematicians ({\it Deutsche  Mathematiker--Vereinigung}), as everywhere in the German 
cultural sphere. In 1862, moreover, a Czech counterpart, the Union of Czech Mathematicians and Physicists 
({\it Jednota \v cesk\'ych matematik\accent23u a fyzik\accent23u}) was founded in Prague by  students. Though this 
{\it Jednota} was initially intended as an organization devoted to the improvement of the students scientific knowledge, 
and teaching and lecturing skills 
without consideration for the language or the ethnic question, the more active Czech members rapidly transformed it into a 
Czech organization, isolating Czech students from the German ones and losing contact with the German teachers. 
The {\it Jednota} became one of the organizations of national consciousness of the Czech intelligentsia.  Though there were members of this society in many places of Bohemia and Moravia, the meetings and lectures were officially held only in Prague till WW1\footnote{See details on the Union of Czech mathematicians and physicists in Seidlerov\'a, I. : Science in a bilingual country in Bohemia in History. Mikul\'a\v s Teich, ed., Cambridge University Press, 1998.}.  

The 1907 of the JDMV issue\footnote{{\it Jahresbericht der Deutschen
Mathematiker--Vereinigung}. 16 1907, pp. 396--397; 18 1909, pp. 104--105; 21
1912, pp. 58--59; 23 1914, pp. 52--53.} mentions 
that in Fall 1905 a mathematico-physical society ({\it Matematisch-Physikalisch Gesellschaft}) was created in Brno and 
that during the academic year 1906-07 a seminar was held with a regular frequency of one  talk every two weeks 
(16 talks that year). Hamel himself read four talks, Waelsch, the other full professor of mathematics in Brno, two, 
and among the other speakers we find the names of M.Ungar, R. von Mises, E.Fischer, E.Fanta\dots 
The full list of the talks between 1906 and 1913 published by the JMDV reveals 63 talks in which other names such as 
Hahn, Haas, Haar, R\"uckle, Ehrenfest, Tietze, Radon appear, some of them several times. An interesting comparison 
can be made with the situation of the mathematical society in Vienna as can be seen in the journal of the DMV. 
Speakers in Brno were often also speakers in the imperial capital the very same year. Tietze, Schrutka, von Mises, Hamel 
came several times to Vienna to speak while they held a position in  Brno, and the mathematical department of the Brno 
German Technical University appeared as a kind of distant suburban university of Vienna. 
One may, however, feel a nuance when comparing the titles of the talks in the two cities. While the Vienna seminar  
concentrated on strictly mathematical aspects, Brno was slightly more oriented towards mechanics and mathematical 
physics (with talks about the Planck's results (Hamel), hydrodynamic (von Mises), electromagnetism (Jaumann), 
gravitation (Jaumann)\dots ) though Radon, Tietze, Fischer and others were also presenting purely mathematical results. 
We can also observe that personal relationships had certainly played a great role in the mathematical seminar in Brno. 
Several mathematicians from the list had met when they had been students in G\"ottingen (Hamel, Fanta and Haar for example). 
Hahn, Tietze and Ehrenfest had formed a small inseparable group of students at the Vienna University. They were in fact four 
in that group to which Gustav Herglotz\footnote{Details about Herglotz and Ehrenfest can be 
found in a recent paper by Huijnen, P. and  Kox, A. J. : Paul Ehrenfest's Rough Road to Leiden : A Physicist's search for a 
Position 1904-1912, Physics in Perspective, 9,2, 186-211, 2007. The same paper presents an interesting description of 
Ehrenfest's travel when he decided to leave St Petersburg in 1912 and to look for a possible permanent position 
somewhere in Europe. His talk in Brno in February 1912 belongs to this period.} also belonged. After Hamel's departure to Aachen 
in 1912, the life of the Gesellschaft, though less active, continued until the beginning of WW1. 

\subsection{Czech institutions}

At the Brno Czech Technical University, several mathematicians were appointed. 
Among the first members of the professors' staff we find  Karel Zahradn\'\i{}k (1848--1916) (who became the first rector) 
and Jan Sobotka (1862-1931).  Zahradn\'\i{}k had studied mathematics and physics at the Prague University, 
then he became a secondary school teacher in Prague and an assistant at the Prague Czech Technical University. 
In 1876 he was appointed Professor of mathematics at the University of Zagreb in Croatia.

As a typical feature, no mathematician or physicist from the Czech
Technical University lectured in the German Mathematical
Society. Until 1860s, the mathematicians 
at universities in Prague and Olomouc were Czechs or Germans born in the Czech
lands, or mathematicians who came to the Czech lands from abroad. They wrote their theses and scientific works in German. 
They were mathematicians
working in Bohemia or Moravia and nobody minded whether they were Czechs or Germans. The situation changed only in the last quarter of the 
19th century, during the time of national revival\footnote{Nov\'y, L. : D\v ejiny exaktn\'\i{}ch v\v ed v \v cesk\'ych
zem\'\i{}ch do konce 19. stolet\'\i{}. Nakladatelstv\'\i{} \v Ceskoslovensk\'e akademie v\v ed,
Praha 1961, p. 221.}. Students from the Prague Technical University called for lectures in Czech. In 1860s the first professors for mathematical lectures
in Czech were appointed and the mathematical community was divided for the first
time. In 1869, the Prague Polytechnicum was divided into two schools and similarly at the beginning of 1880s, the 
venerable Charles University of Prague was split. German mathematicians
from the Czech lands kept tight contacts with mathematicians in Germany
and Austria and were considered as German mathematicians in Europe\footnote{As Seidlerov\'a observes in Seidlerov\'a, I. : Science in a bilingual country in
 {\it Bohemia in History}, the situation continued after the independence of Czechoslovakia and even inside the country. She writes : {\it Even Czech university teachers and researchers
 often had no idea that in their works they actually cited a colleague from the Brno German Technical University}. We besides already described in 
Havlov\'a, V., Mazliak, L. and \v Si\v sma, P. :  Le d\'ebut des relations math\'ematiques franco-tch\'ecoslovaques vu \`a travers la correspondance Fr\'echet-Hostinsk\'y the 
amazing case of F.Urban, a Brno German mathematician who wrote a book in 1923 where he studied random events in chains (Markov chains, in modern terms) at precisely the same time when Hostinsk\'y  became interested in them and who became known to Hostinsk\'y through Fr\'echet.}. 

The Czech mathematicians, though 
more isolated,  were also in contact with European mathematicians. They studied at universities in Germany, France or 
Italy. Their studies generally meant only short, one-year stays after they had  graduated. Notwithstanding this had also 
been the tradition of German mathematicians
from Prague, the mission of Czech mathematicians was new and different as they
had to prepare the first Czech textbooks of mathematics not only for university students, but also for 
secondary school pupils. While Czech mathematicians
wrote textbooks devoted to general mathematics, German
mathematicians could prepare textbooks for specialized mathematical disciplines
and present their scientific results there.  At the end of 19th century, 
the mathematicians, graduating at Czech universities, had a very small chance
to be appointed professors in Austria or Germany. Their biggest hope was 
to be appointed in one of the two Czech universities in Prague. Thus Czech mathematicians mostly worked at secondary 
schools and their conditions for scientific work were very limited.

In 1897,  Jan Sobotka was appointed 
Professor of Descriptive geometry at the Vienna Technical University.  Sobotka had studied in a German real-school 
in Prague. He subsequently specialized in mathematics and descriptive geometry at both Czech universities in Prague. 
He became an assistant at the Czech Technical University where he replaced Til\v{s}er for lectures of 
descriptive geometry. During this period,  he attended seminars of leading German geometers in Zurich and Breslau. 
As he could not obtain a position of teacher even at a secondary school in Prague, he went to Vienna, where he became 
a secondary school teacher at a real-school for a short time. Soon he was appointed extraordinary
Professor of descriptive geometry at the Vienna Technical University. In 1899 he was appointed
Professor in Brno but in 1904 he went to Prague where he received a position of Professor of Geometry at the Prague 
Czech University. The majority of Czech mathematicians of the first half of 20th century  learnt geometry under his 
direction\footnote{Urban, A.,
Van\v cura, Z. : St\'e v\'yro\v c\'\i{} narozen\'\i{} profesora Sobotky. \v Casopis pro p\v estov\'an\'\i{}
matematiky. 87, 1962, pp. 382--386.}.

Another Professor of mathematics at the Brno Czech Technical University
was Anton\'\i{}n Sucharda (1854-1907), a former teacher in a real-school in
Prague. He was also an assistant of descriptive geometry at the Prague Czech
Technical University but he could not obtain a position at universities. He
studied in G\"ottingen, Munich, Paris, and Stra\ss burg. Sucharda worked in Brno
only until 1904 when he fell ill. During the first years of the Brno Czech Technical University, 
V\'{a}clav \v{R}eho\v{r}ovsk\'{y} (1849--1911) came to Brno and was appointed
Professor of mechanics. He had been an assistant of mathematics at the Prague Czech
Technical University and then a teacher at secondary schools.

As we can see, the foundation of the Czech Technical University in Brno resulted in creating four professorships for 
the Czech mathematical community. For two of them, not young secondary school teachers, it was an 
opportunity to  obtain better conditions for their scientific work. For Zahradn\'\i{}k and Sobotka,  it was a 
chance to return to their country. In 1906,  after Sucharda's retirement,  Maty\'{a}\v{s} Lerch (1860--1922) replaced him.

 Lerch  was a brilliant number theorist  who had been appointed Professor at  the University of Fribourg, Switzerland in 1896. 
Lerch had been supported by the French mathematician Charles Hermite to obtain this position because Lerch's chances 
to be appointed in Prague were minimal. Lerch had studied mathematics at the Prague 
University and Technical University where he was asked to become an assistant to Eduard Weyr in 1885 and later to Gabriel 
Bla\v zek. He habilitated at the Prague Technical University in 1886 and during the next ten years published more than 
one hundred mathematical works\footnote{Frank, L. : O \v zivot\v e profesora Maty\'a\v se Lercha. \v Casopis pro p\v estov\'an\'\i{}
matematiky. 78, 1953, pp. 119--137.}.  During his stay in Fribourg Lerch tried, without success, 
to obtain a position at the Prague Czech University after the deaths of
Studni\v cka and Weyr and appears as a good example of a first-rate mathematician specialist who could not  find any position in Czech universities as only two of them existed. In 1906, Lerch was appointed
Professor of Mathematics at the Brno Czech Technical University, but his lectures seemed to have been  more suitable for 
mathematical students than for future engineers. Lerch's assistant and successor Karel \v Cupr wrote in Lerch's obituary\footnote{\v Cupr, K. : Profesor Maty\'a\v s Lerch. \v Casopis pro p\v estov\'an\'\i{}
matematiky a fysiky. 52, 1923, p. 309.} 
that Lerch's lectures in Brno were the same as the ones he had read at the  University in Fribourg. There were even public 
protests of the Technical University students against Lerch. 

Eventually, with the choice of Miloslav Pel\'\i{}\v{s}ek (1855-1940) as Professor of descriptive geometry in 1908, the staff 
situation in mathematics became stabilized\footnote{Fran\v ek, O. : D\v ejiny \v Cesk\'e vysok\'e \v skoly
technick\'e v Brn\v e. Vol. 1, Vysok\'e u\v cen\'\i{} technick\'e v Brn\v e, Brno 1969.
pp. 233--240.}.

The foundation of the second Czech Technical University also brought changes into the organization of the  
{\it Jednota}. The Brno section of {\it Jednota} was officially 
founded in 1913 but the members had started organizing lectures and meetings immediately after the establishment of the Brno 
Czech Technical University. During the years 1901--1911 about 55 lectures were held in Brno.
The auditors were mostly Brno secondary school teachers of mathematics\footnote{Ko\v s\v t\'al, R. : Vznik a
v\'yvoj pobo\v cky J\v CMF v Brn\v e. Jednota \v ceskoslovensk\'ych matematik\accent23u a fyzik\accent23u.
Brno 1967. pp. 18--21.}.

\section{Mathematics in Brno after WW1}

As seen from the previous picture,  though there has been a mathematical life in Brno in the years preceding WW1, it was    
completely split between the two national communities,  as was the case in political questions. We were not able to 
find a single example of mathematical cooperation not only between  a German and a Czech university at an official level (such as exchange of professors, common lectures and so on) but even between two individual members of 
these universities!

The mathematical life in Brno became paralyzed by the beginning of WW1. The activity of German mathematical society 
at the German Technical University stopped immediately. After Hamel's departure from Brno,  
the professorship of theoretical mechanics remained vacant until 1916. Both Professors of mathematics - Lothar Schrutka 
(1881--1945) as well as Heinrich Tietze (1880--1964) - were enlisted in the army. Tietze was an officer in the front, Schrutka taught at a secondary military school 
in Vienna. From the German Technical University in Brno 9 professors, 10 docents, 40 assistants, and 34 others employees 
were enlisted\footnote{Haussner, A. : Geschichte der Deutschen Technischen Hochschule in Br\"unn
1849--1924. In Festschrift der Deutschen Technischen Hochschule in Br\"unn zur
Feier ihres f\"unfundsiebzigj\"ahrigen Bestandes im Mai 1924. Verlag der Deutschen
Technischen Hochschule, Br\"unn 1924. pp. 35--36.}. All mathematical lectures and lessons were read by the Professor of descriptive geometry Emil Waelsch (1863--1927). 
Both his assistants were enlisted and they were replaced  by students. Waelsch's assistants Wilhelm Schmid (1888--1963) 
and Rudolf Kreutzinger (1886--1959) were imprisoned in Russia and returned to Brno only in 1920 and 1921 respectively - which resulted in a difficult situation for descriptive geometry teaching  immediately after the war, when the number of students increased. The number of students decreased from 950 before the war to 100-200 during the war.

Paradoxically, the situation of teaching at the Czech Technical University was better during the first years of the war 
because all the professors of mathematics and descriptive geometry (Zahradn\'\i{}k, Lerch, and
Pel\'\i{}\v{s}ek) were old men, and therefore not enlisted in the army. Zahradn\'\i{}k died in 1916 and his successor, 
Jan Vojt\v{e}ch (1879--1953), was appointed only in 1918. From 177 teachers and employees of the Czech Technical University in Brno,  
64 persons were enlisted\footnote{Fran\v ek, O. : D\v ejiny
\v Cesk\'e vysok\'e \v skoly technick\'e v Brn\v e. Vol. 1, Vysok\'e u\v cen\'\i{} technick\'e v Brn\v e,
Brno 1969. p. 103.}. The number of regular students at the school decreased from 571 in the academic year 1913/14 
to 254 in the academic year 1914/15 and there were approximately 90 students in 1916/17\footnote{Archives of the Czech Technical 
University, Brno. Lists of students of academic years 1913/14 to 1919/20.}. 
The remaining students of these empty years were often younger than before the war and belonged to the classes not yet 
called to the army. In 1913/14, there were 30\% to be under 19 years old, while in 1914/15 the same category represented 
49\% of the students. 
In 1917, a certain number of soldiers were allowed to come back to universities,
and the number of students  increased to 368 :  many of them were students having passed their first year examinations 
in 1914/15 and they were now enrolled in their second year.

The establishment of Czechoslovakia brought a turn in the position of Brno
technical universities. The financial situation of universities had become critical and several buildings of both Brno universities had been 
used as military hospitals. The damage in these buildings made the restart of teaching in 1918 more difficult\footnote{\v Si\v sma, P.
: Matematika na n\v emeck\'e technice v Brn\v e. Prometheus, Praha 2002. pp. 146--148.}.  

 Being at a loss because of 
their new position of political inferiority, the leaders of German universities often made alarming declarations. It is true that as Germans boycotted the Czechoslovak National Assembly of 1919 (which resulted not from elections but from a common agreement of Czech and Slovak leaders)\footnote{The first general election was held in 1920, and 72 German deputies were elected (that is, deputies who belonged to a party whose name included the word `German').}, the Czech majority was at ease to ask for a tight control of the German institutions. While 
the Czechoslovakian parliament  was discussing the organization of 
higher education in the new state, a professor in the Prague German University declared : {\it We are in an appalling situation that a great part of our university would be thrown out onto the street. 
The situation is distressing and is best expressed by the words : homeless, 
without means, without rights}\footnote{{\it Co to znamen\'a? To znamen\'a, \v ze akademick\'y sen\'at n\v emeck\'e university v Praze
neuznal pr\'avn\'\i{} stav tohoto st\'atu, pon\v evad\v z p\v ripou\v st\'\i{} mo\v znost, \v ze by ta
universita mohla b\'yti p\v renesena do st\'atu ciz\'\i{}ho, \v ze tedy n\v emeck\'a universita
v Praze nen\'\i{}, abych tak \v rekl, majetkem nebo statkem \v Ceskoslovensk\'e
republiky}.
Quoted in Parliament Discussion,  27 February 1919 34th Session.}.

In the same session of the Parliament on 27 February 1919,  Srd\'\i{}nko, who was a 
deputy at the Parliament then, contested the honesty of 
such declarations and 
claimed for a substantial reevaluation of the public means to the Czech and German universities. To support his assertion, 
Srd\'\i{}nko mentioned that even before the war,  texts had been published by foreign authors to condemn the disproportion of 
means between the German and Czech universities. He quoted a paper in the 
Revue G\'en\'erale from 1911 where one mentions that {\it a brutal and obvious fact appears from this amount of documents. 
It is the extraordinary disproportion existing between the credits attributed to the German University and those 
attributed to the Czech University, if one takes into account the respective populations.\footnote{Un fait brutal, \'evident, se d\'egage de cette masse des documents.
C'est l'extraordinaire disproportion, qui existe entre les cr\'edits affect\'es
\`a \v l'Universit\'e allemande et ceux qui sont affect\'es \`a \v l'Universit\'e tch\`eque, si
\v l'on tient compte de leur population respective.
}}

The discussions at the Czechoslovakian Parliament in the years 1919-1920 were opportunities to  
present an avalanche of figures aiming at proving that the German universities received satisfactory funding  from the 
Government and had no reasons to complain.  Considering the example of the Brno (German) Technical University, 
Mare\v{s} mentions\footnote{Parliament Discussion, February 27th, 1919 34th Session.} that its budget in the Austrian times 
was 707 thousand crowns, while it was 1753 thousand crowns in 1919. It is of course difficult to  precise the 
exact meaning of such figures as the war had considerably shaken the exchange rates. At the beginning of the war, in order 
to make it popular, the Austro-Hungarian government in Vienna had decided to pay double price for the main articles of 
necessity (grain, cattle, horses\dots ) : an enormous amount of 30 billions Crowns had been printed by the Austro-Hungarian 
bank with forced rate during the war, resulting in a huge inflation\footnote{Ra\v s\'\i{}n, A. : Financial policy of Czechoslovakia 
during the first year of its history, Humphrey Milford, Oxford, 1923,  p. 23.}.
The estimated inflation index rate for the crown in the Czech lands in October 1918 was 1876 (100 in 1914)\footnote{\v Sediv\'y, I. : \v Ce\v si,
\v cesk\'e zem\v e a velk\'a v\'alka 1914--1918, Nakladatelstv\'\i{} Lidov\'e Noviny,  p.
245.}. 

The situation at the German Technical University in Brno was therefore more complicated than at the Czech 
Technical University after WW1. We have seen that before the war, a lot of teachers had come to Brno from Austrian 
universities. These often young men, had not seen real differences between positions in Graz, Brno or Innsbruck, and 
Brno was certainly a more attractive town for them than  Lemberg or Czernovitz. After WW1, the situation changed and the 
Technical University in Brno, though it remained an institution intended for a minority, but now a minority without 
political power, became therefore much less attractive. The number of students was twice greater than before the war at the German Technical University in school year 1920/21, but the future of the German Technical University was obscure.  
 Following the independence declaration of Czechoslovakia on 29 October 1918 
the German deputies from Bohemia and Moravia in Vienna had proposed the future annexation of the regions with German 
settlement to Germany (north Bohemia and Silesia) or Austria (south Bohemia and south Moravia). The Czechoslovakian 
Government solved the problem by sending its army there from the end of November 1918. In March 1919, a violent repression 
against the Sudete Germans temporarily concluded the question and the German inhabitants resigned to belong to the 
Czechoslovakian state\footnote{B\v elina, P., \v Cornej, P. et Pokorn\'y, J. : Histoire des Pays Tch\`eques, Points-Histoire, 
Seuil, 1993.}. 

The existence of two German technical universities for the German minority was nevertheless a political problem 
because the Czechs soon drew attention to the fact that three million Germans had the same number of technical 
universities as nine million Czechs and Slovaks.  The Czechoslovak government chose, however, to avoid a crisis, and 
a statu quo decision was quickly made. In March 1919, the professors of the Brno German Technical University took a vow 
of loyalty to the Czechoslovak Republic\footnote{Archive of the German Technical University in Brno, Moravian Provincial Archive,
B 34, 416.}. 

During this hectic period of difficult political situation, Austrian universities worried about the destiny of their `sister-institutions'  in the Czech lands. On December 14th, 1918 a meeting of Austrian universities was organized in Vienna  to discuss a possibility of a common future.  The representatives of the Prague German University asserted that it would be impossible to continue their 
activity in Prague, and proposed to move to a town in Bohemia where German settlement was in 
majority\footnote{Quoted in Czechoslovak Parliament Discussion, February 27th, 1919 34th Session.}.
At the same time the Association of Austrian German Ingenieurs proposed to transfer the Prague German Technical 
University to \' Usti-nad-Labem in the north of Bohemia and to officially transform Brno (German) Technical University into a branch of Vienna Technical University\footnote{Lidov\'e Noviny, January  29th, 1919.}.
On December 23rd, 1918 the Academic Senate of the Prague German University declared that if the regions with German settlement  obtained their reunion with  Germany or Austria, the Prague University should also transfer.  This declaration infuriated deputy Mare\v s expressed his anger about that point : {\it What does 
that mean? It means that the Academic Senate of the German University in Prague does not acknowledge the legality of this 
country [and acts] as if the German University in Prague were not a property of the Czechoslovakian  
Republic}\footnote{Parliament Discussion, February 27th, 1919 34th Session.}.

The complicated situation at the Brno German Technical University resulted in many changes in the local mathematical 
community.  In 1919, the professor of mathematics Heinrich Tietze left for Erlangen University and Ernst Fanta, 
an actuarial from Vienna, stopped his lectures of actuarial mathematics in Brno where he had come every week until 
that moment. It is possible that other docents from Vienna who had also worked at the Brno Technical University before 
the war, stopped commuting.   The position of Professor of mathematics remained vacant until 1923, when 
Karl Mayr (1884--1940) was appointed. Eight mathematicians applied for this place in 1919 and the professors staff chose 
Johann Radon as the best candidate. Unfortunately, Radon had accepted an offer to become professor at Hamburg University in 
the meantime. The second possible candidate for the professorship in Brno, Roland Weitzenb\"ock, had already been 
appointed at the Prague German Technical University. The negotiation continued in 1921 until the rector of the Brno 
Technical University suggested a professor of mathematics at the Clausthal Mining Academy, Horst von Sanden. Sanden 
rejected the offer, and so did Robert K\"onig (1885--1979) and Georg Prange (1885--1941) in 1922. Eventually Karl Mayr, an  assistant of mathematics 
in Brno before WW1 became a privatdocent at the Vienna Technical University and Tietze's successor. But that lasted for a 
short time. Mayr, dissatisfied with his position of extraordinary professor in Brno, left the town 
for the Graz Technical University in 1924 (though in Graz  he was offered the same status of extraordinary Professor)\footnote{\v Si\v sma, P. : Matematika na n\v emeck\'e
technice v Brn\v e. Prometheus, Praha 2002. pp. 216--219.}.

In 1925, Lothar Schrutka left Brno after having accepted professorship at the Vienna Technical University. The long 
negotiations about Mayr's and Schrutka's successors were again difficult. In 1925, Rudolf Weyrich (1894--1971), 
a student of Breslau University and a privatdocent in Marburg, was appointed extraordinary  Professor and in 1927, 
Lothar Koschmieder (1890--1974) came to Brno as an ordinary Professor\footnote{\v Si\v sma, P. : Matematika na n\v emeck\'e
technice v Brn\v e. Prometheus, Praha 2002. pp. 220--225.}. The discussions about the positions of professors 
at the Brno German Technical University were extremely slow. The 
very bad exchange rate of the Czechoslovak crown promised only very small salaries in comparison with their equivalent in Austria or Germany. In 1927, Josef Krames (1897--1986), a 
Professor of descriptive geometry at the German Technical University in Brno from 1927 to 1929, wrote to the Czechoslovak 
Ministry of Education, that his future salary in Brno would be the same as his salary of an assistant in Vienna\footnote{Czech
National Archive in Prague, Ministry of Education, Josef Krames --- personal
file.}. 
He nevertheless accepted the position, because he rightly hoped that the status of extraordinary Professor in Brno would help him 
to obtain the professorship of descriptive geometry at the Graz Technical University. This happened in 1929.

Another factor delaying the negotiations was the fact that the Czechoslovak government obstinately refused foreign 
candidates to be chosen. The Ministry of Education required exclusive proposals for Czechoslovak experts in Czechoslovak 
universities and only in some exceptional situations (when a native expert could not be found) might a foreign candidate 
be considered. As the German universities of Prague or Brno did not produce a sufficient number of specialists, able to 
compete with Austrian and German ones, they naturally turned towards foreign scholars and the discussions were complicated. 
Let us also add that after WW1 the mathematical life at the German Technical University in Brno obviously could not 
continue the traditions of the German Mathematical Society before the war. Only fundamental mathematical
lectures for engineers remained at the Technical University  and there were no more special lectures read by 
privatdocents of mathematics. The Professors of mathematics at 
the German universities in Czechoslovakia kept their contacts with mathematicians in Austria and Germany and did all they 
could to obtain positions at foreign universities. They regularly participated in the meetings of the 
{\em Deutsche Mathematische Vereinigung} and organized such a meeting in Prague in 1929. While in Prague we know about 
individual contacts between Czech and German mathematicians, in Brno these contacts seem to have been extremely rare.

Rather than a confrontation with the German community,  the Czechoslovak government chose to make a priority of the future development of the Czech Technical
University (as proven by the creation of the new Faculty of 
Architecture in 1919 after many years of requests).  When the peace came back, the number of regular students 
at the Czech Technical University experienced a big jump to more than 900. Moreover,  the Technical universities were opened to women after the war\footnote{Before 1919, women could study at technical universities only as extraordinary students. The admittance of any woman was discussed at the professors staff meetings. Though there were several attempts to open technical universities to women before the Great War, the Austrian Ministry of education remained always opposed to the idea.}. Due to the presence of the older classes, the 
age distribution was slightly different than before the war : in 1914 only 14\% of the students were over 24 years old 
and they were 28\% in 1919.  Many of them were therefore in need of a rapid technical qualification to start working as the technicians of the new state, a situation not particularly favorable to mathematics. And indeed the number of professors of mathematics, descriptive geometry, and mechanics did not change at the Brno technical 
universities after WW1. After Lerch's departure to Masaryk University,  Karel \v{C}upr (1883-1956) was appointed his 
successor in 1923. Descriptive geometry was 
taught by Pel\'\i{}\v sek until 1928. 
Another name to be mentioned,  Jur Hronec (1881--1959) came to Brno in 1924\footnote{Hronec worked at the Technical University till 
1938.  He was the only Slovak professor of mathematics at Czechoslovak universities, illustrating the 
somehow problematic disparity between the two founding  components of the multinational state. In 1939 Hronec was 
appointed Professor and the first rector of the first Slovak Technical University in Ko\v sice. Descriptive geometry was taught by Pel\'\i{}\v sek at Brno Technical
University till 1928 (see Fran\v ek, O. : D\v ejiny \v Cesk\'e vysok\'e \v skoly
technick\'e v Brn\v e. Vol. 1, Vysok\'e u\v cen\'\i{} technick\'e v Brn\v e, Brno 1969. pp.
233--240).}.

\section{The foundation of Masaryk University}

We now arrive at the main event in the academic life of Brno in the immediate afterwar period, the creation of an entirely new university. For the new government, the creation of high education institutions was not only a conclusion to the aforementioned fights for the rights of the Czech citizens. It was now a question of vital necessity to build such institutions in the south of the country because the traditional road to Vienna for Moravian and Slovak students was barred. 
Yet, before the collapse of the monarchy, Czech and Slovak leaders agreed on the necessity of creating two new universities, 
in Brno and in Bratislava.  Not later than at the second session of the 
Czechoslovak Parliament on November 15th, 1918  a group of deputies proposed the 
creation of a university in Brno with three faculties : of philosophy, law and medicine. 
At the beginning of 1919, it was decided to divide the philosophical faculty into two faculties --- philosophical and 
science as had been done in the Prague universities and the institution with four faculties called Masaryk University 
was founded by the 
Law of January 28th, 1919.  Before the beginning of the first academic year,  the first professors of 
each faculty were appointed.  The lectures started at the Faculty of Law (a bit strangely only for the first and third 
years of study --- perhaps because of the unusual age distribution already mentioned in case of the Czech Technical University) 
and at the Faculty of Medicine (the first two years of study) in the academic year 1919/20. The faculties of Science and  
Philosophy started their activities the following year. Among the first professors at the Faculty of Science we find the 
names of Bohuslav Hostinsk\'y (1884--1951), Maty\'a\v s Lerch or Bed\v rich Mack\accent23u (1879--1929).

The new university had to solve many problems. The first task was the designation of the professors board. As may be 
expected, the main source of teachers for Masaryk University was the Prague Czech University and its privatdocents, 
including the mathematicians Bohuslav Hostinsk\'y and Ladislav Seifert (1883--1956). 

Bohuslav Hostinsk\'y was the son of a very famous member 
of the Czech intelligentsia, the musicologist Otakar Hostinsk\'y.  In 1906, he defended a thesis in mathematics under 
the title {\it On Lie spherical Geometry}. For the academic year 1908-1909,   Hostinsk\'y obtained a grant from the 
Ministry of Education for one-year stay in Paris. He followed lectures by Picard, Poincar\'e and Darboux there. His 
Parisian stay was a crucial moment for his scientific evolution and allowed him to prepare his habilitation. Back in 
Prague Hostinsk\'y became again a Gymnasium teacher in 1909-10 and then, from 1910,  in the  {\it Realschule} of 
Prague-Vr\v{s}ovice. This was his first permanent position, an important step for any teacher in the lands with the 
German schooling system.  At the same time, he was finishing his habilitation which he defended on November 16th, 
1911 under the title {\it On Geometric methods in the theory of functions}. In 1912,  Hostinsk\'y was called 
to be privatdocent to the Prague 
University. In parallel with his secondary teaching, he began to give conferences on several themes of higher mathematics 
(theory of analytical functions, differential geometry of curves and surfaces, differential equations, geometric 
applications of differential equations \dots ). For reasons we do not know (most probably health reasons), Hostinsk\'y was not enlisted during the war and remained in Prague. We have narrated how Hostinsk\'y discovered probability theory during the war in another paper\footnote{Havlov\'a, V., Mazliak, L., \v Si\v sma, P. : Le d\'ebut des relations math\'ematiques franco-tch\'ecoslovaques vu \`a travers la correspondance Fr\'echet Hostinsk\'y, Jehps, 1,1. 2005.}.  Some months before his appointment in Brno, during the academic year 1919-1920, Hostinsk\'y taught Volterra's theory on integral equations and their applications\footnote{Ber\'anek, J. : Bohuslav Hostinsk\'y (1884--1951). \v Casopis
pro p\v estov\'an\'\i{} matematiky. 109, 1984, pp. 442--448.}.  We shall see in the next section how Hostinsk\'y played a major r\^ole for developing mathematics in Brno in the inter-wars period.

 Ladislav Seifert was a privatdocent of the Prague University. In the academic year 1907/08 Seifert 
had studied at universities in Stra\ss burg and G\"ottingen. Seifert wrote works devoted to algebraic geometry and differential geometry. His works dealt with surfaces of third order and quadrics of revolution. In differential geometry he studied the properties of some curves and surfaces and he interpreted his results in descriptive geometry. He was therefore an heir of the Czech geometers  of the second half of 19th century and he was outside the main trends of the development of geometry of the 20th century\footnote{On Seifert, see  Hrdli\v ckov\'a, J. : \v Zivot a d\'\i{}lo
Ladislava Seiferta (1883--1956). Thesis. Masaryk University, 2001. About the Czech 
geometric school, see the contribution by Folta in  L'Europe Math\'ematique, C.Goldstein, J.Gray and J.Ritter, Edrs, 
Maison des Sciences de l'Homme, Paris, 1996.}.  In 1920, Seifert habilitated for mathematics at the Prague 
University and at the Technical University for descriptive and
synthetic geometry and it was probably the urgent need for a 
teacher of geometry for the new Masaryk University that allowed Seifert to obtain the position of professor. 

The main brain teaser for the new university was, however, the absence of rooms and
buildings. All the faculties started their teaching in temporary conditions, some rooms
were rented from the Technical University or other organizations. The first
mathematical lectures by Lerch were held in a room  where he had taught his students of the Technical University before. 
A huge plan for developing an academic quarter near the
building of the Technical University was not realized, and in fact only one building of that
quarter was built --- the Faculty of Law\footnote{It is today the building where the archives of the university are located.}. The other faculties were 
located in different parts of the town. The situation did not really change until today. A larger campus is nowadays being constructed.

Masaryk University had to build its own libraries for faculties or departments.  Bohuslav Hostinsk\'y 
was the zealous head of the library commission of the Faculty of Science and managed to obtain the inheritance of several personal libraries to create scientific libraries in Brno. 

Hostinsk\'y spent also a lot of energy founding a journal of the Brno Faculty of Science {\em Spisy 
vyd\'avan\'e p\v r\'\i{}rodov\v edeckou fakultou Masarykovy university}. Hostinsk\'y (as we shall see in the next section) 
managed to obtain many exchange agreements with hundreds of scientific institutions from the whole world. This appeared 
to be a decisive point for the new university as, on the score of the very bad exchange rate of the Czech crown,  the 
foreign journals were unaffordable. Thanks to the exchange, 
it could be declared in 1937 that the department of mathematics bought only six journals,  and obtained one hundred 
through exchanges\footnote{\v Cech, E. : O reorganizaci
na\v s\'\i{} v\v edy. Na\v se v\v eda. 18, 1937.}. Moreover, the creation of this new journal became a good opportunity not only for the faculty members 
but also for all Czechoslovak scientists to disseminate their results. The issues were published as separate numbers. 
In 1925, there had already been about fifty of them. 

In the new Faculty of Science, there were two professorships for mathematics and two for physics. Maty\'a\v s Lerch 
and Ladislav Seifert were appointed the first Professors of mathematics. Seifert was appointed Professor of geometry, Lerch taught mathematical analysis and algebra. 
Bohuslav Hostinsk\'y, a privatdocent at the Prague University, was appointed Professor of theoretical physics, and Bed\v rich Mack\accent23u, who had been extraordinary
Professor of physics at the Czech Technical University in Brno, was appointed Professor of experimental physics.

The creation of Masaryk University enabled Lerch to transfer and to give lectures for 
mathematicians, especially secondary school teachers. In 1920, Lerch shared with Hostinsk\'y the teaching of all mathematical and physical lectures and officially started the 
development of the mathematical department. Nevertheless, as Lerch was ill and tired, his actions were restricted and he 
mostly limited himself to teaching. His main contribution may be to have  noticed the talented student Otakar 
Bor\accent23uvka (1899--1995) whom he managed to appoint the first
assistant of mathematics at Masaryk University. Bor\accent23uvka  then became the leader of Brno mathematics in the 
second half of 20th century.  In August 1922,  Lerch died. Seifert became the head of the mathematical department and Eduard \v Cech (1893--1960) was appointed in 1923 only one year after his habilitation. \v Cech had studied mathematics and descriptive geometry at the Prague University from 1912. He was enlisted during the war, but, for reasons unknown to us,  was not sent to the front. He served as a clerk in the rear for three years, having time for reading mathematical books 
and learning foreign languages  (Italian, German, and Russian). After the war he started to teach in a secondary school 
in Prague. In 1920 he defended a thesis on differential geometry. He spent the academic year 1921/22 in Torino under the 
direction of G. Fubini. Back in Prague, he received the habilitation degree at the Prague University in 1922 and a year 
later he was appointed extraordinary Professor in Brno. He was a geometer but he taught mathematical analysis and algebra. 
In 1928 he was appointed full Professor and became one of the world best specialist in topology during the 1930s\footnote{Kat\v etov, M., Nov\'ak, J.,
\v Svec, A. : Akademik Eduard \v Cech. \v Casopis pro p\v estov\'an\'\i{} matematiky. 85, 1960,
pp. 477--491.}.

\section{An axis Brno-Strasbourg : dreams and reality}

Once Masaryk University settled down, it appeared necessary for it to participate in international scientific programs. 
More than an academic need, this was clearly a strong political symbol of the scientific presence of the new 
Czechoslovakia.

In the mathematical science, a very interesting example of efforts to establish scientific relationships linked with the 
new political map of Europe is given by the contacts between Maurice Fr\'echet (1878-1973)  in Strasbourg and Bohuslav Hostinsk\'y in 
Brno. We shall explain how these effort found a basis in the international situation of the time. This situation was the conjugation of several aspects :  the eagerness of the allies (mostly of the continental ally, 
France\footnote{On the French politics towards central Europe, see Wandycz, P. S. : France and its Eastern Allies 1919-1925, Minneapolis University Press, 1962.}) to consolidate its alliance with the new countries of Central Europe  to thwart a possible future German awakening; the aspirations of Czechoslovakia to be more present on the international stage; the anxiety of Germany to keep tight links with its former natural sphere of influence. 

France probably wanted to make an emblem from its relationship with Czechoslovakia.  A singular fact is that very soon in the French official rhetoric appeared a comparison between Czechoslovakia and Alsace,  two countries presented as having been rescued from the jaws of German imperialism. This idea influenced all the domains of economical 
and cultural life, in particular academic life. In 1919, the new born Czechoslovakia seemed an excellent opportunity for  the program of cultural exchanges that 
France wanted to organize after the war. The new top politicians (beginning by the emblematic Tom\'a\v s Masaryk and his second 
Eduard Bene\v s) had kept tight personal and intellectual contacts with France.  An extremely active propaganda was organized 
by the French authorities to convince the Czech Government and the local administrations (universities, schools, cultural 
associations) of the importance of cultural and educational contacts. Two universities of Brno and Bratislava, newly 
created in 1919, were objects of a special attention. 

After the war, the traditional competition between Germany and France to attract students and scholars from foreign countries had encountered a new development. Victorious France had now a position of conqueror (and the officials often showed an amazing self-confidence on the subject), while Germany was on a defensive line. As soon as in 1918, a German scientist, Rieser, wrote in {\it Akademische Rundschau} that, if the Germans do not provide the necessary efforts to attract foreign students  {\it afterwards, Russians and Japanese will go to French schools which are not worse than the German ones, and will come back home and spread the French spirit} \footnote{Akademische Rundschau, V, 1918, p.322.}.

Almost immediately after the end of the war, the reconquest of the university of Strasbourg, and its reconstruction along 
French standards appeared an urgent task to the French Government.  An interesting sign is given in Lavisse's speech 
(Lavisse was the Director of the Ecole Normale) to the students for the opening of Academic Year 1919-1920 at the Ecole 
Normale Sup\'erieure in Paris : {\it You will be soldiers again, and you will serve with honour in the French intellectual 
army whose center is in Paris, and its advanced guard in Strasbourg}.  

Let us immediately observe that the fact that Strasbourg had strong links with German culture was not  an obvious 
advantage in the eyes of foreign students from the former lands of the central Empires as the French authorities had 
first thought. A French diplomat in Transylvania (a part of the Austro-Hungarian Empire given to Rumania) who also 
tried to attract students to Strasbourg writes to the Dean of Strasbourg university in June 1920 : 

\begin{quote} For many people here, and among the most francophile, the Alsacian remains a hybrid person, German as 
much as French, and who, condemned to live on one of the two sides, finally prefers France where his life is more 
comfortable. [The local academic persons responsible] fear that in Strasbourg one cannot breathe an absolutely pure French air. 
You certainly recognize therein the effect of the Fritz propaganda on the mind of these brave Transylvanians.\footnote{Pour beaucoup de gens ici, et des plus francophiles, l'Alsacien demeure un \^etre hybride, autant allemand que fran\c cais, et qui, r\'eduit \`a vivre dans l'un des deux camps, pr\'ef\`ere en d\'efinitive celui de la France o\`u on lui fait la vie plus douce.J'ai senti [qu'ils craignent] qu'on ne respire pas \`a Strasbourg `un air de France absolument pur'. Vous reconnaissez l\`a, bien certainement, Monsieur le Doyen, les effets de la propagande boche sur l'esprit de ces braves transylvains.} \end{quote}

Nevertheless, the French government lived in hopes of attracting many foreign students from eastern Europe to Strasbourg and wanted the university to become a display of the successes of French science. The following letter from a French deputy to the administrator of Alsace, dated from March 1919, is a good illustration.

\begin{quote}
You know better than anyone the considerable importance that the Germans had given to this university and the coquettishness they showed to make of it one of the most brilliant, if not the most brilliant, of the whole empire. You certainly have read that going away they predicted that in less than 3 years France would have jeopardized their work. How can we face this challenge\footnote{Vous savez mieux que personne l'importance consid\'erable que les allemands avaient donn\'ee \`a cette universit\'e et la coquetterie qu'ils ont mise \`a en faire une des plus brillantes sinon la plus brillante de l'empire. Vous avez certainement vu aussi qu'ils ont pr\'edit en partant qu'en moins de 3 ans la France aurait sabot\'e leur \oe uvre. 
Comment relever ce d\'efi?}?
\end{quote}

Strasbourg, therefore, became a first rank university in France for ten years, and a place of original scientific 
experiments. It is in Strasbourg that  French scientists, who faced the experience of the war organization (especially the 
direction of inventions for national defense where Borel had played a major role) understood the importance of developing statistics. During the imperial period, Stra\ss burg had indeed been a major place of the discipline, with Lexis and Knapp\footnote{On the development of mathematical statistics in France, 
see Catellier, R. and Mazliak, L. : Borel, IHP and the beginning of mathematical statistics in France after WW1. To appear.}. Among the pedagogical initiatives where the teaching of statistics was of prime importance, there was the creation of the Institute for Commercial studies ({\it Institut d'\'etudes commerciales}) where an original teaching was made by the sociologist Maurice Halbwachs and Fr\'echet.

Before being appointed to Strasbourg in 1919, Fr\'echet had been professor at Poitiers University between 1910 and 1914, 
and he had been already famous in the international mathematical community after having defended an outstanding thesis 
on the topology of functional spaces in 1906 in which he offered a theoretical framework for the use of Volterra's 
functions of lines. 
In addition to his mathematical fame, Fr\'echet had another asset : he was a polyglot. A singular fact of  his life is 
the energy he devoted to the promotion of Esperanto : in particular, he wrote several mathematical papers in this language. 
He had an excellent knowledge of English, at a time when this was not usual. During the war, he served as an interpreter 
for the British army. He knew German well, a very useful thing in Strasbourg. 

From the very beginning of his presence in Strasbourg, Fr\'echet had seriously taken the question of 
international relations of the university into consideration. On June 29th, 1919 he wrote the following 
letter to Prague :

\begin{quote}
My dear colleague,

May I ask you to let me know which are the universities that should remain and which should be created on the territory 
of your new state. Moreover maybe one of your students would like to oblige by sending to me the list of professors of 
mathematics of the Czechoslovakian universities, as well as the list of Czechoslovakian journals printing original papers 
in mathematics written by your fellow countrymen mathematicians.  Is any of these journals publishing in French? 

Will you excuse me, my dear colleague, for all these questions. Receive my most respectful regards.

Maurice Fr\'echet 
\end{quote}

It is not clear precisely to whom Fr\'echet had written, but we know that  after several months Sobotka asked K\"ossler to transmit the letter to Hostinsk\'y, who received it on October 19th, 1919. Hostinsk\'y was 
then the secretary of the National Provisory Comittee of Czechoslovakia\footnote{Bru, B. : Souvenirs de Bologne, 
{\it Jour. Soc. Fr. Stat}, 144, 135-226, 2003.}, and this may explain why he was put in charge of answering 
Fr\'echet. As he had stayed in Paris for a while, he was also probably known for his good knowledge of French. 

Hostinsk\'y answered to Fr\'echet on October 19th, 1919 from Prague (it had been a few weeks before he left for  Brno). 
He informed Fr\'echet about the future opening of the universities in Brno and Bratislava, and specified that 
the Brno faculty of law would open soon, whereas it would be so in 1920 for the scientific disciplines. Hostinsk\'y also 
mentioned that the two most important journals, the {\it \v{C}asopis pro p\v{e}stov\'an\'\i \hskip 2pt  matematiky a 
fysiky} (Journal for the cultivation of mathematics and physics) and the {\it V\v{e}stn\'\i k Kr\'alovsk\'e \v{c}esk\'e 
spole\v{c}nosti nauk} (Bulletin of the Royal Czech Science Society) would change soon their language policy and decided to 
increase the presence of French and English to the detriment of German. To conclude his letter, Hostinsk\'y did not forget 
to mention that he had been in Paris during the academic year 1908-09, and that he had studied there with the jewel of 
French mathematics (Darboux, Poincar\'e, Picard, Humbert, Appell, Hadamard, Borel\dots ). He proposed to Fr\'echet to 
become his main contact in Czechoslovakia in case of need. 

 Fr\'echet answered to this letter on November 12th, 1919, lavishing advice on Hostinsk\'y for collecting all French 
abstracts of all Czech publications in a single 
journal. One may feel in Fr\'echet's letter a slight touch of paternalism towards new developing communities. Fr\'echet 
was certainly conscious of the fact, as he cautiously wrote that he had thought {\it interesting to let [Hostinsk\'y] 
know the opinion of a stranger who seeks nothing but good things for Czech scientists and mathematical science}, and that 
collecting these abstracts would show {\it how large the part of Czech science was among what was usually attributed 
to the Germans in Austria.}

In a further letter dated June 1st, 1920, written on a heading paper of the {Organization committee of the sixth 
International congress of Mathematicians},  Fr\'echet enclosed a little brochure called {\it Teaching of mathematics at 
the University of Strasbourg} that had been printed in order to attract students in the Alsace capital, and asked 
Hostinsk\'y to include its publication in  a Czech journal. Fr\'echet wrote that the Strasbourg {\it University needs to 
create new currents towards itself and for a while it will be necessary to make some propaganda}.  Hostinsk\'y answered at 
the end of June. He hoped to meet Fr\'echet for the first time in Strasbourg, as he would be a member of the Czech 
delegation to the International Congress of Mathematicians held in the town in September 1920.  

The Czech delegation was important (11 members, of the total of about 200 persons). In a report made after the congress\footnote{\v Casopis pro
p\v estov\'an\'\i{} matematiky a fysiky. 50, 1921, p. 46--47.}, Bohumil 
Byd\v{z}ovsk\'y (1880--1969) mentioned that \begin{quote} the contact with mathematicians from the Strasbourg University, which is our 
main university partner in the West, was particularly cordial. The interest they showed for our scientific, educational 
and social situation seems to warrant that reciprocal exchanges will continue, obviously for the prosperity of our 
science.\footnote{Zvl\'a\v st\v e srde\v cn\'y byl styk s matematiky \v Strasbursk\'e university, na\v s\'\i{} nejbli\v z\v s\'\i{}
sp\v r\'atelen\'e vysok\'e \v skoly na z\'apad\v e. Z\'ajem, kter\'y jevili o na\v se pom\v ery v\v edeck\'e,
studijn\'\i{} i jin\'e ve\v rejn\'e, zd\'a se b\'yti z\'arukou, \v ze ve vz\'ajemn\'ych styc\'\i{}ch bude
pokra\v cov\'ano, jist\v e ve prosp\v ech na\v s\'\i{} v\v edy.
} \end{quote}

At Strasbourg congress where he indeed met Fr\'echet, Hostinsk\'y read two talks : one on differential geometry, the second 
one on mechanics. Moreover, during the spring of 1920, Hostinsk\'y had sent to Emile Picard the translation of his Czech 
paper, published in 1917 in {\it Rozpravy \v{C}esk\'e Akademie}\footnote{Hostinsk\'y, B. : Nov\'e \v re\v sen\'\i{} Buffonovy \'ulohy o jehle.
Rozpravy \v cesk\'e akademie II. 26, 1917. p. 8.} and devoted to a new solution of Buffon needle problem. 
As soon as he received it (April 18th, 1920), Picard proposed to include the article into the Miscellaneaous section of 
the {\it Bulletin des Sciences Math\'ematiques}. This slightly modified version of the 1917 paper was published at the end of 
1920\footnote{Hostinsk\'y, B. : Sur une nouvelle solution du probl{\`{e}}me de l'aiguille.
Bulletin des Sciences Math\'ematiques. 44, 1920. pp. 126--136.} and Fr\'echet read it carefully, as he himself mentioned in a letter to Hostinsk\'y dated November 7th, 1920 
congratulating the author for having obtained a {\it positive result}. This reading became an opportunity for Fr\'echet 
to write his first paper in probability theory, and at the same time for Hostinsk\'y the first step towards his studies on the 
ergodic principle. 

The beginning of the correspondence with Hostinsk\'y was only one reason for Fr\'echet to become interested in probability and  not the main one. His lectures at Strasbourg Commercial Institute together with Halbwachs clearly played a greater role. However, our guess is that it was the continuation of this correspondence which introduced him to a major probabilistic theme, Markov chains. Hostinsk\'y's decisive luck was indeed a submission he made to the Academy of Sciences in Paris in 
1928. It was a note devoted to an elementary version of the ergodic theorem for a continuous state Markov chain. The 
matter was to meet a spectacular development in the 1930s and Hostinsk\'y preceded there all the future major specialists 
(Kolmogorov in the first place). The importance of the note did not escape Hadamard and marked for the French 
mathematician the only period of his long life devoted to probability theory, a period nicely referred to by Bru as 
Hadamard's {\it ergodic spring}, ended by the International Congress of Bologne (September 1928) where Hadamard  
read a lecture on card shuffling. Between February and June 1928, Hostinsk\'y and Hadamard exchanged a lot of letters, 
published several notes responding to one another and also met during Hadamard's journey to Czechoslovakia in May. From this moment, Hostinsk\'y acquired a real international prestige and in the 1930s his little school in Brno became 
an active research center around markovian phenomena. Until its collapse in the dramatic events of the German annexation 
of 1939 and World War 2, this little school of mathematical physics in Brno was one of the most successful creations of a 
mathematical center directly inherited from the new situation of postwar Europe. 

The letters of the time with  Fr\'echet prove how Hostinsk\'y presented all the developments of the subject to his French colleague. It is besides what Fr\'echet himself wrote in memoriam at Hostinsk\'y's death
\begin{quote}
Among his so varied researches, he had succeeded in introducing me to the theory of probabilities in chain. Hence it is thanks to him that I could write the second volume of my studies on modern probability theory on the subject, and in the book I frequently used his ingenious methods\footnote{Parmi toutes ses recherches si vari\'ees, il avait su m'int\'eresser  \`a la th\'eorie des probabilit\'es en cha\^\i ne. De sorte que c'est  gr\^ace  \`a lui que j'ai \'et\'e amen\'e  \`a \'ecrire sur ce sujet le second livre de mes Recherches sur la th\'eorie moderne des probabilit\'es, ouvrage o\`u j'ai eu  \`a invoquer ses ing\'enieuses m\'ethodes en de nombreux passages\dots  (Fr\'echet to the Rector of Masaryk university on 5 May 1951).}.
\end{quote}

\section*{Conclusion}

The foundation of Czechoslovakia in 1918 appears a good example of an attempt of reorganization of Europe after the end 
of World War I. In places where there was a tense cohabitation of several national communities (as in many parts 
of the collapsed Austro-Hungarian empire), it was necessary to choose a form of organization allowing the coexistence 
of several traditions. This was in particular the case with the organization of educational system. 

When one studies the local case of Brno, the capital of Moravia on the border of Austria, it is 
vital to understand how the difficult contacts between the Czech majority and the large German minority had influenced 
the whole process of edification of education institutions between ca 1880 and 1930. Though the German minority  
lost its domination in 1918, the institutions were still much influenced by the culture that had prevailed before the war, 
though there were several attempts to create a new interest towards the countries of the victorious side. 

Considering the case of mathematics, we tried to expose how the discipline was mainly active in the German Technical 
University before the war with the creation of a local German Mathematical Society when Hamel was given a position of 
full professor of Mechanics, and in the new Masaryk University after the war where we emphasized the important 
role of Hostinsk\'y and his international contacts, especially with France. However, despite political changes 
unsolved contradictions, antagonisms and absence of communication continued to exist between the two communities.

\end{document}